\documentclass[12pt]{article}
\usepackage{amsfonts}
\usepackage{latexsym}
\usepackage{amsmath}
\usepackage{amssymb}
\usepackage{amsthm}
\usepackage{fullpage}
\usepackage{enumerate}
\newtheorem{theorem}{Theorem}[section]
\newtheorem{lemma}[theorem]{Lemma}

\newtheorem{proposition}[theorem]{Proposition}
\newtheorem{corollary}[theorem]{Corollary}
\newtheorem{remark}[theorem]{Remark}
\newtheorem{conjecture}[theorem]{Conjecture}
\theoremstyle{definition}

\newtheorem{thmy}{Theorem}
\newenvironment{oldtheorem}{\stepcounter{thm}\begin{thmy}}{\end{thmy}}
\newtheorem*{note*}{Note}



\begin{document}

\small

\title{\bf More on logarithmic sums of convex bodies}

\medskip

\author {Christos Saroglou}

\date{August, 2014}
\maketitle

\begin{abstract}
\footnotesize We prove that the log-Brunn-Minkowski inequality (log-BMI) for the Lebesgue measure in dimension $n$
would imply the log-BMI and, therefore, the B-conjecture for any log-concave density in dimension $n$. As a consequence,
we prove the log-BMI and the B-conjecture for any log-concave density, in the plane. Moreover, we prove that the log-BMI
reduces to the following: For each dimension $n$, there is a density $f_n$, which satisfies an integrability assumption, so that the
log-BMI holds for parallelepipeds with parallel facets, for the density $f_n$. As byproduct of our methods, we study possible
log-concavity of the function $t\mapsto |(K+_p\cdot e^tL)^{\circ}|$, where $p\geq 1$ and $K$, $L$ are symmetric convex bodies,
which we are able to prove in some instances and as a further application,
we confirm the variance conjecture in a special class of convex bodies. Finally, we establish
a non-trivial dual form of the log-BMI.
\end{abstract}
\section{Introduction}$\textnormal{ }$

Let $K,\ L$ be convex bodies in $\mathbb{R}^n$ (i.e. compact, convex sets, with non-empty interiors), that contains the origin in their interiors.
For $p\geq 1$, the $L^p$-Minkowski-Firey sum $a\cdot K+_pb\cdot L$ of $K$ and $L$ with respect to some positive numbers $a,b$
is defined by its support function
\begin{equation}
h_{a\cdot K+_pb\cdot L}=(ah_K^p+bh_L^p)^{1/p}\ .
\label{Lp-sum,p>1}
\end{equation}
The case $p=1$ corresponds to the classical Minkowski sum $aK+bL=\{ax+by|x\in K, y\in L\ \}$. In the pioneer
work of B\"{o}r\"{o}czky, Lutwak, Yang and Zhang \cite{BLYZ}, the $L^p$-convex combination of $K$ and $L$, with respect to some
$\lambda\in(0,1)$, for all $p\geq 0$ is defined:
$$\lambda \cdot K+_p(1-\lambda)\cdot L=\{x\in \mathbb{R}^n\ |\ x\cdot u\leq [\lambda h^p_K(u)+(1-\lambda)h^p_L(u)]^{1/p},
\textnormal{ for all }u\in S^{n-1}\} \ \ \ \textnormal{and}
$$
$$\lambda \cdot K+_0(1-\lambda)\cdot L=\{x\in
\mathbb{R}^n\ |\ x\cdot u\leq h^{\lambda}_K(u)h^{1-\lambda}_L(u), \textnormal{ for all }u\in S^{n-1}\} \ .
$$
Note that if $0\leq p<1$, $\lambda \cdot K+_p(1-\lambda)\cdot L$ cannot be defined by (\ref{Lp-sum,p>1}), simply because the resulting
function is not always convex. Nevertheless, the two definitions coincide for $p\geq 1$.
Let us state the fundamental $L^p$-Brunn-Minkowski inequality  (see e.g \cite{Ga2}, \cite{Sch}, \cite{Ga1}),
for $p\geq 1$ in its dimension-free form:
$$|\lambda \cdot K+_p(1-\lambda)\cdot L|\geq |K|^{\lambda}|L|^{1-\lambda}\ ,$$
where $|\cdot |=|\cdot |_n$ is the $n$-dimensional Lebesgue measure.

Although for $p\geq 1$, the $L^p$-Brunn-Minkowski theory has been considerably developed in the previous years
(see e.g. \cite{Lut1} \cite{Lut2} \cite{LYZ1}
\cite{LYZ2} \cite{LYZ3}), much less seem to be known for $0\leq p<1$.
The following is conjectured in \cite{BLYZ} (without the equality cases):
\begin{conjecture}\label{Conjecture-log-bm-lebesque}(The logarithmic-Brunn-Minkowski inequality)
Let $K$, $L$ be symmetric convex bodies in $\mathbb{R}^n$ and $\lambda\in(0,1)$. Then,
$$
|\lambda\cdot K+_o(1-\lambda)\cdot L|\geq |K|^{\lambda}|L|^{1-\lambda}\ ,
$$
with equality in the following case: Whenever $K=K_1\times\dots\times K_m$, for some convex sets $K_1,\dots ,K_m$, that cannot be written
as cartesian products of lower dimensional sets, then there exist positive numbers $c_1,\dots,c_m$, such that
$L=c_1K_1\times\dots\times c_mK_m \ .$
\end{conjecture}
The conjecture can easily be seen to be wrong for general convex bodies, even for $n=1$. Note, also, that for $0\leq p\leq q$,
$\lambda \cdot K+_p(1-\lambda)L\subseteq \lambda \cdot K+_q(1-\lambda)L$, thus the log-Brunn-Minkowski inequality (if true) implies the
$L^p$-Brunn-Minkowski inequality for all $p>0$.

The cone-volume measure of $K$ is defined as: $S_0(K,\cdot)=h_KS(K,\cdot)$, where $S(K,\cdot)$ is the surface area measure of $K$, viewed
as a measure on the sphere (see e.g. \cite{Sch}). In \cite{BLYZ2}, a necessary and sufficient condition was discovered (the planar
case was treated by Stancu \cite{St1} \cite{St2}; see also \cite{BLYZ3} for other applications of the cone-volume measure and \cite{Co-Fr}
for a possible functional generalization of the classical Minkowski problem). A confirmation of Conjecture \ref{Conjecture-log-bm-lebesque}
would answer the following open problem: When do two symmetric convex bodies $K,L$ have proportional cone-volume measures? If
Conjecture \ref{Conjecture-log-bm-lebesque} was proven to be true, the pairs $(K,L)$ would be exactly the ones for which equality
holds in the log-Brunn-Minkowski inequality. The planar case was settled in \cite{BLYZ}:
\begin{oldtheorem}\label{planar-log-BM}\cite{BLYZ} Conjecture \ref{Conjecture-log-bm-lebesque} is true in dimension two.
\end{oldtheorem}
It was shown in \cite{Sa1} that
Conjecture \ref{Conjecture-log-bm-lebesque} holds true for pairs of unconditional bodies with respect to the same orthonormal basis.
Actually, the proof (based on a result from \cite{CFM}) shows that this result (as for the inequality) remains true if we replace the
Lebesgue measure with any unconditional log-concave measure in $\mathbb{R}^n$.
Recall that a measure $\mu$ is called log-concave if for all convex bodies $K$, $L$, it satisfies the
Brunn-Minkowski inequality:\[\mu\big(\lambda K+(1-\lambda)L\big)\geq \mu(K)^{\lambda}\mu(L)^{1-\lambda}\ , \ \lambda\in(0,1)\ .\]
By a result of C. Borell \cite{borel}, the absolutely continuous log-concave measures in $\mathbb{R}^n$ are exactly the ones having
log-concave densities, i.e. their logarithms are concave functions.
It is reasonable to conjecture the following:
\begin{conjecture}\label{Conjecture-log-bm-general}
Let $\mu$ be an even log-concave measure in $\mathbb{R}^n$, $K$, $L$ be symmetric convex bodies and $\lambda\in (0,1)$. Then,
$\mu\big(\lambda\cdot K+_o(1-\lambda)\cdot L\big)\geq \mu(K)^{\lambda}\mu(L)^{1-\lambda}\ .$
\end{conjecture}
Conjecture \ref{Conjecture-log-bm-general} is closely connected (actually implies; see Corollary \ref{log-BM implies general B-conjecture})
with the so called B-conjecture.
\begin{conjecture}\label{conjecture-b-conjecture}(B-conjecture)
Let $\mu$ be an even log-concave measure in $\mathbb{R}^n$ and $K$ be a symmetric convex body. Then, the function
$\mathbb{R}\ni t\mapsto \mu(e^tK)$
is log-concave.
\end{conjecture}
This was conjectured in \cite{CFM} \cite{Kl3}. This was previously conjectured by  Banaszczyk in \cite{La} for the standard
Gaussian measure $\gamma_n$ (i.e. the measure that has density $e^{-\|x\|_2^2/2}$)
and it was confirmed by Cordero-Erasquin, Fradelizi and Maurey in \cite{CFM},
a fact known as the B-theorem (see also \cite{La-Ol} for an application).
More generally it was shown that
\begin{oldtheorem}\label{B-theorem-Gaussian}(The B-Theorem for the Gaussian measure \cite{CFM})
Let $A$ be a diagonal $n\times n$ matrix and $K$ be a symmetric convex body. Then, the function
$\mathbb{R}\ni t\mapsto \gamma_n(e^{At}K)\ $is log-concave.
In particular, the standard Gaussian measure satisfies the B-conjecture.
\end{oldtheorem}
Moreover, the following fact, also from \cite{CFM}, will be used:
\begin{oldtheorem}\label{B-theorem-unconditional}\cite{CFM}
Let $A$ be a diagonal $n\times n$ matrix, $\mu$ be an unconditional log-concave measure and $K$ be an unconditional convex body. Then, the function
$\mathbb{R}\ni t\mapsto \mu(e^{At}K)$
is log-concave.
\end{oldtheorem}
A connection between the log-Brunn-Minkowski inequality for the Lebesgue measure (Conjecture \ref{Conjecture-log-bm-lebesque}
and the B-conjecture for uniform measures of symmetric convex bodies (i.e. measures of the form $|K\cap\cdot|$, where $K$ is a symmetric
convex body) was established in \cite{Sa1}. Namely, it was proven that (i) the log-Brunn-Minkowski inequality for the Lebesgue measure in
dimension $n$
implies the B-conjecture for uniform measures in dimension $n$. Thus, by Theorem \ref{planar-log-BM}, the B-conjecture for uniform measures
in the plane follows (this fact was proven independently in \cite{Li}). (ii)
The log-Brunn-Minkowski inequality for the Lebesgue measure holds in any dimension if and only if, in any dimension,
the function $|(e^{At}C_n)\cap K|$ is
log-concave in $t$ for any symmetric convex body $K$ and for any diagonal matrix $A$.
Here $C_n$ denotes the cube $[-1,1]^n$.

Our first goal is to continue the ideas from \cite{Sa1} and
extend the formentioned results even further. Let us briefly describe our main results towards this direction.
In Section 3 we prove (see Theorem \ref{log-BM for Lebesgue implies general log-BM}) that actually the log-Brunn-Minkowski inequality
for the Lebesgue density
implies the log-Brunn-Minkowski inequality for any log-concave density and, therefore, the B-conjecture in full generality.
Thus, again by Theorem \ref{planar-log-BM}, we establish (see Corollary \ref{log-BM-planar-general})
Conjectures \ref{Conjecture-log-bm-general} and
\ref{conjecture-b-conjecture} in the plane.

On the other hand, in Section 5 we modify the proof of fact (ii) mentioned earlier to prove that actually
in order to confirm
Conjecture \ref{Conjecture-log-bm-general}, one needs for any dimension $n$ to find a density $f_n$ which satisfies a mild integrability
assumption and the function \[t\mapsto\int_{e^{At}C_n}f_n(Tx)dx\]
is log-concave for any choice of the diagonal matrix $A$ and for any invertible linear map $T$. The reader should focus in the case of the
Gaussian density; see Remark \ref{remark-section-3}.

Our second goal is to study log-concavity and log-convexity properties for dual bodies. 
In Section 4, as byproduct of our method from Section 3, we show (see Proposition \ref{Section2-proposition1}) that the B-conjecture
for uniform measures or for measures with densities of the form $e^{-\|x\|_K^p}$, $p\geq 1$ would imply the log-concavity of the function
$$t\mapsto|(K+_p\cdot e^tL)^{\circ}|\ ,$$where $M^{\circ}$ stands for the dual body of $M$. Using the cases where the B-Conjecture is known
to hold, we establish this log-concavity property in some special cases (see Corollary \ref{Section2-corollary2}).
As a further application, in Theorem \ref{Section2-Thm-Poincare}
we confirm the variance conjecture (see Section 4 for more information) in a special class of convex bodies.

Finally, in Section 6, we establish the $L^0$-analogue
of Firey's dual Brunn-Minkowski inequality \cite{Fi1} (and its extension to other quermassintegrals):
\begin{equation*}
 \big|\big(\lambda K+(1-\lambda)L\big)^{\circ} \big|\leq |K^{\circ}|^{\lambda}|L^{\circ}|^{1-\lambda}\ .
\label{Firey-intro}
\end{equation*}
Note that the $L^0$-version is clearly a stronger inequality. Also, since no explicit formula is 
valid for the support function of the logarithmic sum, no classical arithmetic
inequalities (such as H\"{o}lder) can be used directly towards the proof.
Therefore, our inequality is a non-trivial extension of Firey's result.

\section{Preliminaries}
\hspace*{1.6 em}Let us state some results that will be needed subsequently. We refer to \cite{Sch} \cite{Ga2} for more information.

Let $K$ be a convex body that contains 0 in its interior. The polar body of $K$ is defined as:
$$K^{\circ}=\big\{x\in\mathbb{R}^n\ |\ x\cdot y\leq 1,\ \forall y\in K\big\}\ .$$ Then,
$K^{\circ}$ is also a convex body that contains 0 in its interior and  $\big(K^{\circ}\big)^{\circ}=K$.

The $i$-th quermassintegral $W_i(K)$ of $K$ is defined by the Steiner formula
\[|K+tB_2^n|=\sum_{i=0}^n\bordermatrix{& \cr &n \cr &i}W_i(K)t^i \  ,\ \ t>0\ ,\ \ i=0,\dots ,n\ \ ,\]
where $B_2^n$ is the Euclidean unit ball. Note that $W_0(K)$ is the volume of $K$ and $W_1(K)$, $W_{n-1}(K)$, $W_n(K)$ are proportional to the surface area,
the mean width and the Euler characteristic respectively. Moreover, the functional $W_i$ is $(n-i)$-homogeneous, that is $W_i(tK)=t^{n-i}W_i(K)$.
A useful formula for the quermassintegrals of $K$ is the Kubota recursion formula:
\[W_i(K)=\int_{{\cal{G}}_{n,n-i}}|K|H|_{n-i}\ dH\ , \ i=1,\dots,n-1 \ .\]
Here, $\int_{{\cal{G}}_{n,n-i}}\cdot\ dH$ denotes the integral of a function defined on the Grassmannian ${\cal{G}}_{n,n-i}$,
with respect to the Haar-measure
on ${\cal{G}}_{n,n-i}$ and $K|H$ is the orthogonal projection of $K$ onto the subspace $H$.

For the rest of this section, $K$ will be denoting a symmetric convex body (i.e. $K=-K$).
The norm $\|\cdot \|_K$ of $K$ is the unique norm in $\mathbb{R}^n$, such that $K=\{x\ |\ \|x\|_K\leq 1\}$. Recall that every norm in
$\mathbb{R}^n$ is the norm of a unique symmetric convex body.

The support function of $K$ is defined by $h_K(x)=\max_{y\in K}(x\cdot y)\ , \ x\in\mathbb{R}^n\ .$
There is a duality relation between the norm and the support function of $K$: $h_K=\|\cdot\|_{K^{\circ}} \ .$

The inradius and the outradius of $K$ are defined as:
\[\textnormal{inradius}(K)=\min_{x\in \partial K}\|x\|_2 \ ,\ \textnormal{outradius}(K)=\max_{x\in K}\|x\|_2 \ .\]

Let $f:\mathbb{R}^n\rightarrow \mathbb{R}_+$ be a homogeneous of degree $p$ function. Then, by integration in polar coordinates, we have:
\[\int_Kf(x)dx=\frac{1}{n+p}\int_{S^{n-1}}f(x)\|x\|_K^{-n-p}dx \ ,\]
where $S^{n-1}=\{x\in\mathbb{R}^n\ | \ \|x\|_2=1\}$, the unit sphere in $\mathbb{R}^n$.

The Pr\`{e}kopa-Leindler inequality is probably the most famous functional generalization of the Brunn-Minkowski inequality
which states that whenever
$f,g,h$ are non-negative measurable functions with the property that for some
$\lambda\in (0,1)$, $h(\lambda x+(1-\lambda)y)\geq f^{\lambda}(x)g^{1-\lambda}(y)$,
for all $x,y\in\mathbb{R}^n$, then
\begin{equation}
\int_{\mathbb{R}^n}h(x)dx\geq \Big[\int_{\mathbb{R}^n}f(x)dx\Big]^{\lambda}\Big[\int_{\mathbb{R}^n}g(x)dx\Big]^{1-\lambda}\ .
\label{eq-pl}
\end{equation}
We will need the 1-dimensional Pr\`{e}kopa-Leindler inequality in the following form, proven in \cite{U}:
\begin{oldtheorem}
\label{oldthm-PL}
(Multiplicative version of the Pr\`{e}kopa-Leindler inequality)\\
Let $\lambda\in (0,1)$, $f,g,h:\mathbb{R}_+\rightarrow \mathbb{R}_+$ be non-negative measurable functions, such that
$h(x^{\lambda}y^{1-\lambda})\geq f^{\lambda}(x)g^{1-\lambda}(y)$, for all $x,y>0$. Then, (\ref{eq-pl}) holds.
\end{oldtheorem}
The proof follows by applying the Pr\`{e}kopa-Leindler inequality to the functions $\overline{f}(x)=e^xf(e^x)$, $\overline{g}(x)=e^xg(e^x)$,
$\overline{h}(x)=e^xh(e^x)$ and the change of variables $y=e^x$.

It is well-known
that there exists a unique-up to isometry-volume preserving linear map $T$ such that the
quantity $$L_{TK}^2:=\frac{1}{|K|^{\frac{n+2}{n}}}\int_{TK}(x\cdot y)^2dx$$ is constant as a function of $y\in S^{n-1}$. Then,
$TK$ is said to be isotropic and the number $L_{TK}$ is called the
isotropic constant of $K$ (see \cite{MP} for basic results on this
concept). It is true that
\[L_K^2=\frac{1}{n|K|^{\frac{n+2}{n}}}\min_{T\in SL_n}\int_{TK}\|x\|_2^2dx\ .\]It has been conjectured that the isotropic constants of symmetric convex bodies are
bounded from above by an absolute constant; this problem is known as the slicing problem. The isotropic constant is known to be bounded form below
by an absolute constant (see again \cite{MP}).
The best estimate up to date for the upper bound is of the order $n^{1/4}$, due to Klartag \cite{Kl1} after improving the previous estimate
$Cn^{1/4}\log n$ by Bourgain \cite{Bou}.

Let $H$ be a $k$-dimensional subspace of $\mathbb{R}^n$. Define the Schwartz-symmetrization $S_H(K)$ of $K$ with respect to
$H$ as the set that is
constructed by replacing every cross-section, orthogonal to $H$, of $K$ with a Euclidean ball of the same $(n-k)$-dimensional volume.
It is an easy application of
the Brunn-Minkowski inequality that $S_H(K)$ is also a convex body.
If $H=\mathbb{R}u$, for some unit vector $u$, we abbreviate $S_u(K)=S_H(K)$.
Notice, furthermore, that
$h_K(u)=h_{S_{u}(K)}(u).$

Let $f:\mathbb{R}^n\rightarrow \mathbb{R}$ be a function. The epigraph of $f$ is defined as
\[Epi(f):=\big\{ (x,t)\ | \ x\in\mathbb{R}^n,\ t\geq f(x)\big\}\subseteq \mathbb{R}^{n+1}\ .\]
It is true that $f$ is convex if and only if its epigraph is a convex set. Moreover, $Epi(f)$ characterizes $f$.
For $u\in S^{n-1}$, define the
Schwartz-symmetrization $S_u(f)$ with respect to $u$, as the function with
\[Epi\big(S_u(f)\big)=S_H\big(Epi(f)\big) \ ,\]
where $H$ is the subspace spanned by $u$ and an orthogonal to $\mathbb{R}^n\equiv$Domain of $f$, vector of $\mathbb{R}^{n+1}$.
By the previous discussion,
if $f$ is convex, then $S_u(f)$ is convex as well.
\section{On the log-Brunn-Minkowski inequality for general log-concave measures}
The main result of this section is the following:
\begin{theorem}\label{log-BM for Lebesgue implies general log-BM}
Assume that the log-Brunn-Minkowksi inequality holds in dimension $n$ for the Lebesgue measure. Then, the log-Brunn-Minkowksi inequality
holds in dimension $n$
for any even log-concave density.
\end{theorem}
\begin{corollary}
\label{log-BM implies general B-conjecture}
Assume that the log-Brunn-Minkowksi inequality holds in dimension $n$ for the Lebesgue measure. Then, the B-conjecture holds in dimension $n$,
for any even
log-concave density.
\end{corollary}
Proof. Let $\mu$ be an even log-concave measure, $K$ be a symmetric convex body, $\lambda\in (0,1)$ and $s,t\in \mathbb{R}$. Then, by Theorem
\ref{log-BM for Lebesgue implies general log-BM}, the log-Brunn-Minkowski inequality holds for the measure $\mu$, therefore
$$\mu\big(e^{\lambda s+(1-\lambda)t}K\big)=\mu\big(\lambda\cdot (e^sK)+_0(1-\lambda)\cdot (e^tK)\big)\geq \mu(e^sK)^{\lambda}\mu(e^tK)^{1-\lambda}$$
and the assertion follows. $\Box$
\\
\\
Combining Theorem \ref{log-BM for Lebesgue implies general log-BM}, Corollary \ref{log-BM implies general B-conjecture} and Theorem \ref{planar-log-BM},
we immediately obtain:
\begin{corollary}\label{log-BM-planar-general}
The log-Brunn-Minkowski inequality and the B-conjecture hold in the plane, for any even log-concave density.
\end{corollary}
Corollary \ref{log-BM-planar-general} (in particular, the B-Theorem in the plane) together with \cite[Proposition 3.1]{Mar} (see also \cite{GZ}) immediately imply the following:
\begin{corollary}\label{Gaussian-BM}
Let $\mu$ be an even log-concave measure in the plane, $M$ be a symmetric convex body in the plane and $\lambda \in (0,1)$. Then for every $K,L \in \{\alpha M ; \alpha \geq 0 \}$, one has
$$ \mu(\lambda K + (1-\lambda) L)^{\frac{1}{2}} \geq \lambda \mu(K)^{\frac{1}{2}} + (1-\lambda) \mu(L)^{\frac{1}{2}}. $$
\end{corollary}
For the proof of Theorem \ref{log-BM for Lebesgue implies general log-BM}, some geometric lemmas are required.
\begin{lemma}
\label{Section 1, first lemma}
Let $\varphi:\mathbb{R}^n\rightarrow \mathbb{R}\cup\{\infty\}$ be an even convex function and $t\in\mathbb{R}$,
so that the sets $\{\varphi\leq t\}$ and $\{\varphi=\varphi(0)\}$ are convex bodies. Then, there exists $b>0$, depending only on $\varphi,t$, such that for all
$u\in S^{n-1}$, $r,s\in\mathbb{R}$ with $\varphi(0)<r<s\leq t $, the following inequality holds:
$$s-r\leq b\big(h_{\{\varphi\leq s\}}(u)-h_{\{\varphi\leq r\}}(u)\big)\ .$$
\end{lemma}

Proof. Note that the restriction of $\varphi$ into the set $\{\varphi\leq t\}$ is Lipschitz with some constant $A>0$. Let $r<s$, $u\in S^{n-1}$ and
$x\in \mathbb{R}^n$, such that
$\varphi(x)=r$ and $x\cdot u=h_{\{\varphi\leq r\}}(u)$. Then, there exists a $\xi>1$, such that for $x':=\xi x$, $\varphi(x')=s$. Then,
$x'\cdot u\leq h_{\{\varphi\leq s\}}(u)$. Thus,
\begin{eqnarray*}
h_{\{\varphi\leq s\}}(u)-h_{\{\varphi\leq r\}}(u)\geq x'\cdot u-x\cdot u
&=&\|x'-x\|_2\frac{x'-x}{\|x'-x\|_2}\cdot u\\
&=&\|x'-x\|_2\frac{x}{\|x\|_2}\cdot u\\
&=&\|x'-x\|_2\frac{h_{\{\varphi\leq r\}}(u)}{\|x\|_2} \ .
\end{eqnarray*}
We have $\varphi(x)=r\leq t$, therefore $\|x\|_2\leq \textnormal{outradius}(\{\varphi\leq t\})$. Also,
$h_{\{\varphi\leq r\}}(u)\geq h_{\{\varphi=\varphi(0)\}}(u)
\geq \textnormal{inradius}(\{\varphi=\varphi(0)\})$. Hence,
\begin{eqnarray*}
h_{\{\varphi\leq s\}}(u)-h_{\{\varphi\leq r\}}(u)&\geq&\frac{\varphi(x')-\varphi(x)}{A}
\frac{\textnormal{inradius}(\{\varphi=\varphi(0)\})}{\textnormal{outradius}(\{\varphi\leq t\})}\\
&=&\frac{s-r}{A}
\frac{\textnormal{inradius}(\{\varphi=\varphi(0)\})}{\textnormal{outradius}(\{\varphi\leq t\})} \ . \ \Box
 \end{eqnarray*}
\begin{lemma}\label{Section 1, second lemma}
Let $\varphi:\mathbb{R}^n\rightarrow \mathbb{R}\cup\{\infty\}$ be a function 
and $t$ be a real number, both satisfying the conditions of Lemma \ref{Section 1, first lemma}. There
exists $c>|\varphi(0)|$, so that if we set
\begin{equation}
\overline{\varphi}(x)=\begin{cases}
\varphi(x)+c & \text{,  }\ \ \varphi(x)\leq t \\
\infty & \text{,  }\ \ \varphi(x)>t \ \ \ ,
\end{cases}
\label{definition of kladiki}
\end{equation}
then for any $r,s\geq 0$,
$$\{\overline{\varphi}\leq r^{\lambda}s^{1-\lambda}\}\supseteq
\lambda\cdot\{\overline{\varphi}\leq r\}+_0(1-\lambda)\cdot\{\overline{\varphi}\leq s\}\ .$$
\end{lemma}
Proof. Let $b$ be the constant from Lemma \ref{Section 1, first lemma}. Set
$$c:=b\cdot\textnormal{outradius}(\{\varphi\leq t\})+|\varphi(0)|$$
and define $\overline{\varphi}$ by (\ref{definition of kladiki}). Clearly, $\overline{\varphi}\geq 0$.
It suffices to prove that for every $r,s>0$,
\begin{equation}
h_{\{\overline{\varphi}\leq r^{\lambda}s^{1-\lambda}\}}(u)\geq h_{\{\overline{\varphi}\leq r\}}^{\lambda}(u)
h_{\{\overline{\varphi}\leq s\}}^{1-\lambda}(u) \ .
\label{Section 1, second lemma, inequality}
\end{equation}
First assume that $\overline{\varphi}(0)\leq r,s\leq t +c$. Fix $u\in S^{n-1}$ and consider the Schwartz symmetrization $S_u(\overline{\varphi})$
of $\overline{\varphi}$.
Note that $$h_{\{\overline{\varphi}\leq p\}}(u)=h_{\{S_u(\overline{\varphi})\leq p\}}(u) \ ,$$
for every $p\geq \overline{\varphi}(0)$. Moreover, since the body $\{S_u(\overline{\varphi})\leq p\}$ is unconditional with respect to an orthonormal basis
that contains $u$, one can easily see that
$$S_u(\overline{\varphi})\Big(h_{\{S_u(\overline{\varphi})\leq p\}}(u)u\Big)=p=S_u(\overline{\varphi})\Big(h_{\{\overline{\varphi}\leq p\}}(u)u\Big)\ .$$
For $t+c\geq p\geq \overline{\varphi}(0)  $, set
$$f(p):=S_u(\overline{\varphi})(pu)\ .$$
Then, $f$ is a convex and strictly increasing function and also,
\begin{equation}
h_{\{\overline{\varphi}\leq p\}}(u)=f^{-1}(p)\ .
\label{Section 1, second lemma, ineq-2}
\end{equation}
We will show that for $t+c\geq s>r\geq \overline{\varphi}(0)$,
$$f^{-1}(r^{\lambda}s^{1-\lambda})\geq f^{-1}(r)^{\lambda}f^{-1}(s)^{1-\lambda}\ .$$
Consider the line through the points $(f^{-1}(r),r)$ and $(f^{-1}(s),s)$ and suppose that this is defined by the equation
$x_2=c_1x_1+d_1$, for points $(x_1,x_2)$ of the plane. Since $f$ is strictly increasing, it is clear that $c_1>0$. We claim that $d_1\geq 0$.
Indeed, by (\ref{Section 1, second lemma, ineq-2}), we have:
\begin{eqnarray*}
d_1=s-\frac{s-r}{f^{-1}(s)-f^{-1}(r)}f^{-1}(s)&=&s-\frac{s-r}{h_{\{\overline{\varphi}\leq s\}}(u)-h_{\{\overline{\varphi}\leq s\}}(u)}
h_{\{\overline{\varphi}\leq s\}}(u)\\
&\geq&s-bh_{\{\overline{\varphi}\leq s\}}(u)\\
&\geq&s-bh_{\{\overline{\varphi}\leq t+c\}}(u)\\
&\geq &s-b\cdot \textnormal{outradius}(\{\varphi\leq t\})\\
&\geq&\overline{\varphi}(0)-b\cdot \textnormal{outradius}(\{\varphi\leq t\})\\
&=&\varphi(0)+|\varphi(0)|\geq 0\ .
\end{eqnarray*}
Now, the convexity of $f$ implies that
$$r^{\lambda}s^{1-\lambda}\leq c_1f^{-1}(r^{\lambda}s^{1-\lambda})+d_1\ .$$It follows that
$$f^{-1}(r^{\lambda}s^{1-\lambda})\geq (r/c_1)^{\lambda}(s/c_1)^{1-\lambda}-d_1/c_1\geq 
\Big(\frac{r}{c_1}-\frac{d_1}{c_1}\Big)^{\lambda}\Big(\frac{s}{c_1}-\frac{d_1}{c_1}\Big)^{1-\lambda}
=f^{-1}(r)^{\lambda}f^{-1}(s)^{1-\lambda}\ .$$
This proves (\ref{Section 1, second lemma, inequality}) in the case where $t+c\geq s,r\geq \overline{\varphi}(0)$.
If $s<\overline{\varphi}(0)$ (or $r<\overline{\varphi}(0)$), then $\{\overline{\varphi}\leq s\}=\emptyset$ and (\ref{Section 1, second lemma, inequality})
holds trivially. On the other hand, if $s>t+c$, then $\{\overline{\varphi}\leq s\}=\{\overline{\varphi}\leq t+c\}$ and
$$\{\overline{\varphi}\leq r^{\lambda}s^{1-\lambda}\}\supseteq\{\overline{\varphi}\leq r^{\lambda}(t+c)^{1-\lambda}\}\ .$$
Thus, if $r\leq t+c$, we fall in the previous cases, otherwise
$\{\overline{\varphi}\leq r\}\supseteq\{\overline{\varphi}\leq t+c\}$ and (\ref{Section 1, second lemma, inequality}) is again trivial. $\Box$

\begin{lemma}\label{Section 1, thrid lemma}
Let $\lambda\in (0,1)$, $a_1,a_2>0$, $\mu$ be a measure, $\varphi:\mathbb{R}^n\rightarrow  \mathbb{R}\cup\{\infty\}$ be a non-negative
even convex function and $K$, $L$ be symmetric convex bodies.
Assume that for all $r_1,r_2>0$,
$$\mu\Big(\big[\lambda\cdot K+_0(1-\lambda)\cdot L\big]\cap \{\varphi\leq r_1^{\lambda}r_2^{1-\lambda}\}\Big)\geq
\mu\Big(K\cap \{\varphi\leq r_1\}\Big)^{\lambda}\mu\Big(L\cap \{\varphi\leq r_2\}\Big)^{1-\lambda}\ .$$
Then,
$$\int_{\lambda\cdot K+_0(1-\lambda)\cdot L}e^{-a_1^{\lambda}a_2^{1-\lambda}\varphi(x)}d\mu(x)\geq
\bigg[\int_{K}e^{-a_1\varphi(x)}d\mu(x)\bigg]^{\lambda}\bigg[\int_{L}e^{-a_2\varphi(x)}d\mu(x)\bigg]^{1-\lambda}\ .$$
\end{lemma}
Proof. By the Fubini Theorem we have:
\begin{eqnarray*}
 \int_{\lambda\cdot K+_0(1-\lambda)\cdot L}e^{-a_1^{\lambda}a_2^{1-\lambda}\varphi(x)}d\mu(x)&=&\int_0^{\infty}
 \mu\Big(\big[\lambda\cdot K+_0(1-\lambda)\cdot L\big]\cap \{ e^{-a_1^{\lambda}a_2^{1-\lambda}\varphi}\geq s\}\Big)ds\\
 &=&\int_{0}^{1}  \mu\Big(\big[\lambda\cdot K+_0(1-\lambda)\cdot L\big]\cap \{ a_1^{\lambda}a_2^{1-\lambda}\varphi\leq -\log s\}\Big)ds\\
 &=&\int_{0}^{\infty}  \mu\Big(\big[\lambda\cdot K+_0(1-\lambda)\cdot L\big]\cap \{ a_1^{\lambda}a_2^{1-\lambda}\varphi\leq r\}\Big)e^{-r}dr\ .\\
\end{eqnarray*}
Set $h(r)=\mu\big(\big[\lambda\cdot K+_0(1-\lambda)\cdot L\big]\cap \{ a_1^{\lambda}a_2^{1-\lambda}\varphi\leq r\}\big)e^{-r}$,
$f(r)=\mu\big(K\cap \{ a_1\varphi\leq r\}\big)e^{-r}$ and $g(r)=\mu\big(L\cap \{ a_2\varphi\leq r\}\big)e^{-r}$. We will make use of the multiplicative
form of the Pr\`{e}kopa-Leindler inequality. If $r_1,r_2>0$, using our assumption, we have:
\begin{eqnarray*}
h(r_1^{\lambda}r_2^{1-\lambda})&=& \mu\Big(\big[\lambda\cdot K+_0(1-\lambda)\cdot L\big]\cap
\{ \varphi\leq (a_1^{-1}r_1)^{\lambda}(a_2^{-1}r_2)^{1-\lambda}\}\Big)e^{-r_1^{\lambda}r_2^{1-\lambda}}\\
&\geq&\mu\Big(K\cap \{\varphi\leq a_1^{-1}r_1\}\Big)^{\lambda}\mu\Big(L\cap \{\varphi\leq a_2^{-1}r_2\}\Big)^{1-\lambda}e^{-r_1^{\lambda}r_2^{1-\lambda}}\\
&\geq&\mu\Big(K\cap \{\varphi\leq a_1^{-1}r_1\}\Big)^{\lambda}\mu\Big(L\cap \{\varphi\leq a_2^{-1}r_2\}\Big)^{1-\lambda}
e^{-[\lambda r_1+(1-\lambda)r_2]}\\
&=&\bigg[\mu\Big(K\cap \{a_1\varphi\leq r_1\}\Big)e^{-r_1}\bigg]^{\lambda}\bigg[\mu\Big(L\cap \{a_2\varphi\leq r_2\}\Big)e^{-r_2}\bigg]^{1-\lambda}\\
&=&f(r_1)^{\lambda}g(r_2)^{1-\lambda} \ .
\end{eqnarray*}
Thus, by Theorem \ref{oldthm-PL}, we have
\begin{eqnarray*}
\int_{\lambda\cdot K+_0(1-\lambda)\cdot L}e^{-a_1^{\lambda}a_2^{1-\lambda}\varphi(x)}d\mu(x)
&=&\int_0^{\infty}h(r)dr\\
&\geq&\bigg[\int_0^{\infty}f(r)dr\bigg]^{\lambda}\bigg[\int_0^{\infty}g(r)dr\bigg]^{1-\lambda}\\
&\geq&\bigg[\int_{K}e^{-a_1\varphi(x)}d\mu(x)\bigg]^{\lambda}\bigg[\int_{L}e^{-a_2\varphi(x)}d\mu(x)\bigg]^{1-\lambda}\ .\ \Box
\end{eqnarray*}
Proof of Theorem \ref{log-BM for Lebesgue implies general log-BM}\\
It is clearly sufficient (by approximation) to prove the log-Brunn-Minkowski inequality for densities of the form 
$e^{-\varphi_{s,t}}$, where $\varphi$ is any even convex function
defined in $\mathbb{R}^n$, $t>s>\varphi(0)$ and $\varphi_{s,t}$ is given by
$$\varphi_{s,t}(x)=\begin{cases}
\varphi(x) & \text{,  }\ \ s\leq\varphi(x)\leq t \\
s & \text{,  }\ \ \varphi(x)< s \\
\infty & \text{,  }\ \ \varphi(x)>t \ \ \ \ \  .
\end{cases}
$$
Note that $\varphi_{s,t}$ satisfies the conditions of Lemma \ref{Section 1, first lemma}, for all choices of $s$ and $t$. Fix $s,t$ and let $c>0$ and
$\overline{\varphi}$ be as in Lemma \ref{Section 1, second lemma}, i.e. $\overline{\varphi}(x)=\varphi_{s,t}(x)+c$. Then,
$$\int_{\lambda\cdot K+_0(1-\lambda)\cdot L}e^{-\varphi_{s,t}}dx=e^c\int_{\lambda\cdot K+_0(1-\lambda)\cdot L}e^{-\overline{\varphi}}dx \ .$$
Therefore, we need to prove the log-Brunn-Minkowski inequality for the even log-concave density $e^{-\overline{\varphi}}$. Note that $\overline{\varphi}$
is non-negative. We need to show that the assumption of Lemma \ref{Section 1, thrid lemma} is satisfied. Let $r_1,r_2>0$. By Lemma \ref{Section 1, second lemma},
we have:
\begin{eqnarray*}
\Big[\lambda\cdot K+_0(1-\lambda)\cdot L\Big]\cap \big\{\overline{\varphi}\leq r_1^{\lambda}r_2^{1-\lambda}\big\}
&\supseteq&\Big[\lambda\cdot K+_0(1-\lambda)\cdot L\Big]\cap\Big[\lambda\cdot\{\overline{\varphi}\leq r_1\}+_0
(1-\lambda)\cdot\{\overline{\varphi}\leq r_2\}\Big]\\
&\supseteq&\lambda\cdot\Big(K\cap\{\overline{\varphi}\leq r_1\}\Big)+_0
(1-\lambda)\cdot\Big(L\cap\{\overline{\varphi}\leq r_2\}\Big) \ .
\end{eqnarray*}
Since we assumed that the log-Brunn-Minkowski inequality holds for the Lebesgue measure, the assertion follows
by taking volumes in the previous inclusion
and by Lemma \ref{Section 1, thrid lemma} (used with $\mu=|\cdot |$ and $a_1=a_2=1$). $\Box$

\section{Log-concavity properties for dual bodies }
\begin{lemma}\label{Section2-lemma1}
Let $M$ be a symmetric convex body and $\mu$ be a measure, such that the function
$$\mathbb{R}\ni t \mapsto\mu(e^tM)$$
is log-concave. Then, for $p\geq 1$, $a\in\mathbb{R}$, the function
$$\mathbb{R}\ni t\mapsto\int_{\mathbb{R}^n}e^{-e^{at}\|x\|^p_M}d\mu(x)$$
is also log-concave.
\end{lemma}
Proof. This is a special case of Lemma \ref{Section 1, thrid lemma}. 
Indeed, set $K=\mathbb{R}^n=L$ and $\varphi(x)=\|x\|_M^p$. For $\lambda\in(0,1)$,
$r_1,r_2>0$, we have:
\begin{eqnarray*}
\mu\Big(\big[\lambda\cdot K+_0(1-\lambda)\cdot L\big]\cap \{\varphi\leq r_1^{\lambda}r_2^{1-\lambda}\}\Big)
&=&\mu\big(\{\|\cdot\|_M^p\leq r_1^{\lambda}r_2^{1-\lambda}\}\big)\\
&=&\mu\big(r_1^{\lambda/p}r_2^{(1-\lambda)/p}M\big)\\
&\geq&\mu\big(r_1^{1/p}M\big)^{\lambda}\mu\big(r_2^{1/p}M\big)^{1-\lambda}\\
&=&\mu\big(\{\|\cdot\|_M\leq r_1^{1/p}\}\big)^{\lambda}\mu\big(\{\|\cdot\|_M\leq r_2^{1/p}\}\big)^{1-\lambda}\\
&=&\mu\big(K\cap\{\varphi\leq r_1\}\big)^{\lambda}\mu\big(L\cap\{\varphi\leq r_2\}\big)^{1-\lambda} \ .
\end{eqnarray*}
Therefore, by Lemma \ref{Section 1, thrid lemma}, if $a_1,a_2>0$,
$$\int_{\mathbb{R}^n}e^{-a_1^{\lambda}a_2^{1-\lambda}\|x\|_M^p}d\mu(x)\geq
\bigg[\int_{\mathbb{R}^n}e^{-a_1\|x\|_M^p}d\mu(x)\bigg]^{\lambda}\bigg[\int_{\mathbb{R}^n}e^{-a_2\|x\|_M^p}d\mu(x)\bigg]^{1-\lambda}\ ,$$
proving our claim. $\Box$
\\
\\
The following is well known. We include its simple proof for the sake of completeness.
\begin{lemma}\label{Section2-lemma2}
Let $M$ be a convex body that contains the origin in its interior. For $p>0$, there exists a constant $c_{n,p}>0$, that depends only on $n$ and $p$,
such that $$|M|=c_{n,p}\int_{\mathbb{R}^n}e^{-\|x\|_{M}^p}dx \ .$$
\end{lemma}
Proof. Write
\begin{eqnarray*}
\int_{\mathbb{R}^n}e^{-\|x\|_{M}^p}dx &=&\int_0^{\infty}\big|\{e^{-\|\cdot\|_{M}^p}>s\}\big|ds\\
&=&\int_0^1\big|\{\|\cdot\|_{M}\leq (-\log s)^{1/p}\}\big|ds\\
&=&\int_0^1\big|(-\log s)^{1/p}M\big|ds\\
&=&\int_0^1(-\log s)^{n/p}|M|ds
=c_{n,p}^{-1}|M|\ .\ \Box
\end{eqnarray*}
\begin{proposition}\label{Section2-proposition1}
Let $p\geq 1$ and $L$ be a symmetric convex body which has one of the following two properties:
\begin{enumerate}[i)]
\item The measure with density $e^{-\|\cdot\|_L^p}$ satisfies the B-Theorem.
\item The uniform measure of $L$ satisfies the B-Theorem.
\end{enumerate}
Then, for any symmetric convex body $K$, the function
$$\mathbb{R}\ni t\mapsto |(K^{\circ}+_pe^t\cdot L^{\circ})^{\circ}|$$
is log-concave.
\end{proposition}
Proof. Suppose that (i) holds. Using Lemma \ref{Section2-lemma1} with $d\mu=e^{-\|x\|_L^p}dx$, we obtain that the function
$$\mathbb{R}\ni t\mapsto \int_{\mathbb{R}^n}e^{-e^{-t}\|x\|_K^p}d\mu(x)=
\int_{\mathbb{R}^n}e^{-\|x\|_L^p-e^{-t}\|x\|_K^p}dx=:\phi(t)\ .$$
is log-concave.
Note that $\big(\|\cdot\|_L^p+e^{-t}\|\cdot\|_K^p\big)^{1/p}$ is the support function (=dual norm) of the convex body $L^{\circ}+_pe^{-t}\cdot K^{\circ}$. Therefore, by Lemma \ref{Section2-lemma2},
$\phi(t)=c_{n,p}^{-1}|(K^{\circ}+_pe^tL^{\circ})^{\circ}|$ and the function
$$|(K^{\circ}+_p\cdot e^tL^{\circ})^{\circ}|=e^{nt/p}|(L^{\circ}+_pe^{-t}\cdot K^{\circ})^{\circ}|$$ is log-concave.\\
Assume now that (ii) holds. Use again Lemma \ref{Section2-lemma1} with $d\mu=\mathbf{1}_K(x)dx$ to get that the function
$$\mathbb{R}\ni t\mapsto \int_Ke^{-e^t\|x\|_L^p}dx$$ is log-concave. Write
\begin{eqnarray*}
\int_Ke^{-e^t\|x\|_L^p}dx&=&\int_Ke^{-\|e^{t/p}x\|_L^p}dx\\
&=&
e^{\frac{-nt}{p}}\int_{e^{-t/p}K}e^{\|x\|_L^p}dx\\
&=&:e^{-\frac{nt}{p}}\mu'(e^{-\frac{t}{p}}K)\ ,
\end{eqnarray*}
where $\mu'$ is the measure with density $e^{-\|x\|_L^p}$. Since the function $e^{-\frac{nt}{p}}$ is log-affine, it follows that the assumption
of Lemma \ref{Section2-lemma1} holds with $\mu'$ instead of $\mu$, thus the function
$$\mathbb{R}\ni t \mapsto \int_{\mathbb{R}^n}e^{-e^t\|x\|^p_K}d\mu'(x)=\phi(t) $$
is log-concave, where $\phi(t)$ was defined previously and was proven to be proportional to $|(K^{\circ}+_pe^t\cdot L^{\circ})^{\circ}|$. This proves our claim. $\Box$
\begin{remark}\label{Section2-remark-after-proposition1}
Proposition \ref{Section2-proposition1} asserts that the B-conjecture for uniform measures implies the log-concavity of the function $t\mapsto
|(K^{\circ}+_pe^t\cdot L^{\circ})^{\circ}|$.
The opposite is also true, since the limiting case $p=\infty$ is just the B-conjecture for uniform measures.
\end{remark}
Next, let us confirm the B-conjecture for uniform measures in its most simple case: The case of the symmetric strips.
\begin{theorem}\label{Secion2-Bthm-for-strips}
Let $u\in S^{n-1}$, $a>0$. Set $E=\{x\in \mathbb{R}^n\ |\ |x\cdot u|\leq a\}$. Then, for every symmetric convex body $K$,
the function $\mathbb{R}\ni t\mapsto |E\cap e^tK|$ is log-concave.
\end{theorem}
Proof. We need to prove that for every $\lambda\in (0,1)$, $t_1,t_2\in\mathbb{R}$, then:
\begin{equation}
|E\cap e^{\lambda t_1+(1-\lambda)t_2}|\geq |E\cap e^{t_1}K|^{\lambda}|E\cap e^{t_2}K|^{1-\lambda} \ . \label{Section-2-B-thm-for-strips-ineq}
\end{equation}
One can easily verify that for each $b>0$,
$$S_u(E\cap bK)=E\cap bS_uK \ , $$
hence nothing changes in (\ref{Section-2-B-thm-for-strips-ineq}) in terms of volumes if we replace $K$ with the Schwartz symmetrization $S_uK$. But then,
$E$, $S_uK$ are unconditional with respect to some (any) orthonormal basis that contains $u$. Now, Theorem \ref{B-theorem-unconditional} proves our claim. $\Box$
\\
\\
It follows immediately by Proposition \ref{Section2-proposition1}, Theorem \ref{Secion2-Bthm-for-strips}, Theorem \ref{planar-log-BM}
and Theorems \ref{B-theorem-Gaussian},
\ref{B-theorem-unconditional}
that:
\begin{corollary}
\label{Section2-corollary2}
Let $K,L$ be symmetric convex bodies, $p\geq 1$ and $u$ be a unit vector. The
function $\mathbb{R}\ni t\mapsto\big|(K^{\circ}+_pe^t\cdot L)^{\circ}\big|$ is log-concave (at least) in the following cases:
\begin{enumerate}[i)]
\item $p=2$ and $L=B_2^n$.
\item $L$ is an origin symmetric line segment.
\item $K$ and $L$ are unconditional, with respect to the same orthonormal basis.
\item $K$ and $L$ are planar.
\end{enumerate}
\end{corollary}

The variance conjecture \cite{A-B-P} \cite{Bo-Ko} states that if $X$ is a random vector with log concave probability density $f$,
whose barycenter is at the origin and
its covariance matrix is the identity (i.e. $X$ is isotropic), then the variance of $\|X\|_2^2$ satisfies
$$Var(\|X\|_2^2)\leq Cn \ ,$$
where $C>0$ is an absolute constant. The variance conjecture plays a central role in modern convex geometry. Surprisingly, it implies
other major conjectures (see \cite{El} \cite{El-Kl}), such as the slicing problem and the KLS conjecture \cite{KLS} up to a logarithmic factor.
The best general known estimate up to date is of order $n^{5/3}$, due to O. Guedon and E. Milman \cite{Gu-Mi}
(see also \cite{Fr-Gu-Pa}). It has been confirmed for random vectors with unconditional log-concave densities \cite{Kl2}
(see also \cite{Fl}, \cite{Co-Go}). We refer to \cite{Gu} for more information and references.

We would like to restrict our attention
in the class of symmetric convex bodies, i.e. the density $f$ is the indicator function of a symmetric convex body.
In this case the variance conjecture
becomes: Let $K$ be a symmetric isotropic convex body. Then,
$$\sigma^2(K):=\frac{|K|\displaystyle\int_K\|x\|_2^4dx-\bigg[\displaystyle\int_K\|x\|_2^2dx\bigg]^2}{\dfrac{1}{n}\bigg[\displaystyle\int_K\|x\|^2dx\bigg]^2}
\leq C\ .$$

\begin{lemma}\label{Section2-lemma-poincare}
Let $K$, $L$ be symmetric convex bodies, $p\geq 1$, $a>0$. Set $T:=(K^{\circ}+_pa\cdot L^{\circ})^{\circ}$. If the function
$$\mathbb{R}\ni t\mapsto \big|\big(K^{\circ}+_pe^ta\cdot L^{\circ}\big)^{\circ}\big|$$
is log-concave, then
$$|T|\int_T\|x\|_L^{2p}dx-\bigg[\int_T\|x\|_L^{p}dx\bigg]^2\leq \frac{p}{a(n+p)}|T|\int_T\|x\|_L^{p}dx \ .$$
\end{lemma}
Proof. Set $f(t)=\big|\big(K^{\circ}+_pe^ta\cdot L^{\circ}\big)^{\circ}\big|$. Then, $f(0)=|T|$ and $f$ is log-concave. 
Integrating in polar coordinates we obtain:
$$f(t)=\frac{1}{n}\int_{S^{n-1}}\Big(\|x\|^p_K+e^ta\|x\|^p_L\Big)^{-n/p}dx \ ,$$
thus
$$f'(t)=\frac{1}{n}\int_{S^{n-1}}\frac{-n}{p}e^ta\|x\|^p_L\Big(\|x\|^p_K+e^ta\|x\|_L^p\Big)^{-(n+p)/p}dx \ .$$
So,
$$f'(0)=-\frac{1}{p}\int_{S^{n-1}}a\|x\|^p_L\|x\|_T^{-(n+p)}dx=-a\frac{n+p}{p}\int_T\|x\|_L^pdx\ .$$
Also,
$$f''(t)=f'(t)+\frac{1}{n}\int_{S^{n-1}}\frac{n(n+p)}{p^2}e^{2t}a^2\|x\|_L^{2p}\Big(\|x\|_K^p+e^ta\|x\|_L^p\Big)^{-(n+2p)/p}dx\ .$$
Therefore,
\begin{eqnarray*}
f''(0)&=&-\frac{a(n+p)}{p}\int_T\|x\|_L^pdx+\frac{a^2(n+p)(n+2p)}{p^2}\int_{S^{n-1}}\frac{1}{n+2p}\|x\|_L^{2p}\|x\|_T^{-(n+2p)}dx\\
&=&-\frac{a(n+p)}{p}\int_T\|x\|_L^pdx+\frac{a^2(n+p)(n+2p)}{p^2}\int_T\|x\|_L^{2p}dx\\
&\geq&-\frac{a(n+p)}{p}\int_T\|x\|_L^pdx+\Big[\frac{a(n+p)}{p}\Big]^2\int_T\|x\|_L^{2p}dx\ .
\end{eqnarray*}
Now, the log-concavity of $f$ implies $f''(0)f(0)\leq [f'(0)]^2$
and the assertion follows. $\Box$
\\
\\
For $a>0$, define the class of convex bodies ${\cal{C}}_a$ as follows:
$${\cal{C}}_a=\Big\{\big(K^{\circ}+_2a\cdot B_2^n\big)^{\circ}\
|\ K\textnormal{ is a symmetric convex body}, \ \big(K^{\circ}+_2a\cdot B_2^n\big)^{\circ}\textnormal{ is isotropic}\ \Big\} \ .$$
Combining Lemma \ref{Section2-lemma-poincare} with Corollary \ref{Section2-corollary2}, we immediately obtain:
\begin{theorem}\label{Section2-Thm-Poincare}
Let $T\in{\cal{C}}_a$, with $|T|=1$, for some $a>0$. Then,
$$\sigma^2(T)\leq \frac{2}{a(n+2)L^2_T}\ .$$In particular, if $a>c/n$, for some absolute constant $c>0$, then $T$ satisfies the variance conjecture.
\end{theorem}
Before ending this section, we would like to give an alternative description of the class ${\cal{C}}_a$.
\begin{lemma}\label{Section2-Description-of-C_a}
Let $K$ be a symmetric convex body and $a$ be a positive number. Then, $(K+_2a\cdot B^n_2)^{\circ}$ is isotropic if and only if
\begin{equation}
\big|(K+_2a\cdot B_2^n)^{\circ}\big|=\max_{T\in SL_n} \big|(TK+_2a\cdot B_2^n)^{\circ}\big|.
\label{Section2-pre-last-equation}
\end{equation}

\end{lemma}
Proof. 
It is easy to check that the quantity $\big|(TK+_2a\cdot B_2^n)^{\circ}\big|$ indeed attains a maximum, among 
$T\in SL(n)$.  
Let $v\in S^{n-1}$, $t\in\mathbb{R}$, $|t|<1$. Define the linear map
$$T_t(x)=\bigg(\frac{1}{1+t}\bigg)^n\big(x+t(x\cdot v)v\big) \ .$$
Then, $T_t\in SL_n$. Using polar coordinates, one may compute:
\begin{eqnarray}
\frac{\partial}{\partial t}\bigg|_{t=0}\Big|\big(T^{-1}_tK+_2a\cdot B_2^n\big)^{\circ}\Big|& =&\nonumber
\frac{\partial}{\partial t}\bigg|_{t=0}\Big|\big(K+_2a\cdot T_tB_2^n\big)^{\circ}\Big|\nonumber\\
&=&\frac{\partial}{\partial t}\bigg|_{t=0}\frac{1}{n}\int_{S^{n-1}}\Big(h_K(x)^2+a\|T_tx\|_2^2\Big)^{-n/2}dx\nonumber\\
&=&\frac{\partial}{\partial t}\bigg|_{t=0}\frac{1}{n}\int_{S^{n-1}}\Big(h_K(x)^2+a(1+t)^{-2/n}\|x+t(x\cdot v)\|_2^2\Big)^{-n/2}dx\nonumber\\
&=&\frac{1}{n}\int_{S^{n-1}}\frac{n}{2}\Big((2a/n)\|x\|_2^2-2a(x\cdot v)^2\Big)\Big(h_K(x)^2+a\|x\|_2^2\Big)^{-(n+2)/2}dx\nonumber\\
&=&(n+2)a\bigg[-\int_{(K+_2aB_2^n)^{\circ}}(x\cdot v)^2dx+\frac{1}{n}\int_{(K+_2aB_2^n)^{\circ}}\|x\|_2^2dx\bigg] \ .\label{Section-2-last-equation}
\end{eqnarray}
Therefore, if (\ref{Section2-pre-last-equation}) holds, then the derivative at $t=0$ of the volume of $\big(T^{-1}_tK+_2a\cdot B_2^n\big)^{\circ}$
equals zero, so by (\ref{Section-2-last-equation}),
$$\int_{(K+_2aB_2^n)^{\circ}}(x\cdot v)^2dx=\frac{1}{n}\int_{(K+_2aB_2^n)^{\circ}}\|x\|_2^2dx \ .$$
Since this is true for all $v\in S^{n-1}$, it follows that $(K+_2a\cdot B^n_2)^{\circ}$ is isotropic. On the other hand, if
$T_0$ is a critical point of the function $SL_n\ni T\mapsto \big|(TK+_2a\cdot B_2^n)^{\circ}\big|$, we have proved that
$(T_0K+_2a\cdot B^n_2)^{\circ}$ is isotropic. By the uniqueness-up to isometry-of the isotropic position, it follows that this critical point is unique,
thus if $(K+_2a\cdot B^n_2)^{\circ}$ is isotropic, then (\ref{Section2-pre-last-equation}) holds. $\Box$

\section{Reduction to the log-BM inequality for coordinate parallelepipeds}

\begin{theorem}\label{thm-reduction}\textnormal{ }
\begin{enumerate}[i)]
\item Assume that for all $n\in\mathbb{N}$, there exists an even function $f_n:\mathbb{R}^n\rightarrow \mathbb{R}$, whose restriction in any subspace
of $\mathbb{R}^n$ is integrable,
with the following property:
For all $T\in GL_n$ and for all diagonal $n\times n$-matrices $A$, the function
$$\mathbb{R}\ni t\mapsto \int_{e^{At}C_n}f_n(Tx)dx$$
is log-concave. Then, the log-Brunn-Minkowski inequality holds for all even log-concave
densities $g:\mathbb{R}^n\rightarrow \mathbb{R}$, for all $n\in\mathbb{N}$.
\item Assume that for all $n\in\mathbb{N}$, there exists an even function $f_n:\mathbb{R}^n\rightarrow \mathbb{R}$, whose restriction in any subspace
of $\mathbb{R}^n$ is integrable, with the following property:
The log-Brunn-Minkowski inequality holds for the density $f_n$, for any two parallelepipeds with parallel facets. Then,
the log-Brunn-Minkowski inequality holds for all even log-concave
densities $g:\mathbb{R}^n\rightarrow \mathbb{R}$, for all $n\in\mathbb{N}$.
\item Fix $n\in\mathbb{N}$. If the log-Brunn-Minkowski inequality holds for some density $f:\mathbb{R}^n\rightarrow \mathbb{R}$,
then for all $T\in GL_n$ and for all
diagonal $n\times n$-matrices $A$, the function
$$\mathbb{R}\ni t\mapsto \int_{e^{At}C_n}f(Tx)dx$$
is log-concave.
\end{enumerate}
\end{theorem}
Proof. Let us first prove (iii). It is easily verified (see \cite{Sa1} \cite{BLYZ}) that if $s,t\in\mathbb{R}$ and $A$ is
a diagonal $n\times n$-matrix, then for
$\lambda\in (0,1)$, $$e^{[\lambda s+(1-\lambda)t]A}C_n=\lambda\cdot \big(e^{sA}C_n\big)+_0(1-\lambda)\cdot\big(e^{tA}C_n\big) \ .$$
Therefore, for $T\in GL_n$,
\begin{equation}
\lambda\cdot \big(Te^{sA}C_n\big)+_0(1-\lambda)\cdot\big(Te^{tA}C_n\big)=
T\big[\lambda\cdot \big(e^{sA}C_n\big)+_0(1-\lambda)\cdot\big(e^{tA}C_n\big)\big]=Te^{[\lambda s+(1-\lambda)t]A}C_n\ .
\label{Section-3-equation-T^-1}
\end{equation}
Thus,
\begin{eqnarray*}
\int_{e^{[\lambda s+(1-\lambda)t]A}C_n}f(Tx)dx
&=&|\det T| \int_{Te^{[\lambda s+(1-\lambda)t]A}C_n}f(x)dx\\
&=&|\det T| \int_{\lambda\cdot \big(Te^{sA}C_n\big)+_0(1-\lambda)\cdot\big(Te^{tA}C_n\big)}f(x)dx\\
&\geq&|\det T|\bigg[\int_{Te^{sA}C_n}f(x)dx\bigg]^{\lambda}\bigg[\int_{Te^{tA}C_n}f(x)dx\bigg]^{1-\lambda}\\
&=&\bigg[\int_{e^{sA}C_n}f(Tx)dx\bigg]^{\lambda}\bigg[\int_{e^{tA}C_n}f(Tx)dx\bigg]^{1-\lambda}\ .
\end{eqnarray*}

Assertion (ii) is just a reformulation of (i). Indeed, one can check that if $P_1$, $P_2$ are two parallelepipeds with parallel facets, then
there exist $s_1,s_2\in\mathbb{R}$, a diagonal matrix $A$ and a $GL(n)$-map $T$, such that $P_i=Te^{sA}C_n$, $i=1,2$. Thus, by
(\ref{Section-3-equation-T^-1}),
$$\int_{\lambda\cdot P_1+_0(1-\lambda)\cdot P_2}f_n(x)dx=|\det T|\int_{e^{[\lambda s+(1-\lambda)t]A}C_n}f_n(Tx)dx\ .$$

It remains to prove (i). Let $K$, $L$ be symmetric convex bodies in $\mathbb{R}^n$ and $\lambda\in (0,1)$. As in \cite[Theorem 1.5]{Sa1}, consider the following
discretized version of the logarithmic sum of $K$ and $L$: Let $v_1,\dots ,v_m$ be unit vectors in $\mathbb{R}^n$, $m\geq n$. Set
$$R_{\lambda}:=\big\{x\in\mathbb{R}^n\ \big|\ |x\cdot v_i|\leq r_i^{\lambda}s_i^{1-\lambda} \ , \ i=1,\dots ,m\big\} \ ,$$
where $r_i=h_K(v_i)$, $s_i=h_L(v_i)$, $i=1,\dots ,m $.
We will prove that under the assumption of (i),
\begin{equation}
|R_{\lambda}|\geq |R_0|^{1-\lambda}|R_1|^{\lambda} \ .
\label{Eq1-section3}
\end{equation}
Since $R_{\lambda}$ can be chosen to be arbitrarily close to $\lambda\cdot K+_0(1-\lambda)\cdot L$, as $m\rightarrow \infty$,
if (\ref{Eq1-section3}) is proved for any choice of the $v_i$'s, $r_i$'s, $s_i$'s, then the log-Brunn-Minkowski inequality will be established for the
Lebesgue measure. But then, by Theorem \ref{log-BM for Lebesgue implies general log-BM},
the log-Brunn-Minkowski inequality for any log-concave measure will follow. Therefore, it suffices to prove
(\ref{Eq1-section3}) for any choice of $m$, $r_i>0$, $s_i>0$, $v_i\in S^{n-1}$, $i=1,\dots, m$. As in \cite{Sa1}, write
$$|R_{\lambda}|=\int_{\mathbb{R}^n}\prod_{i=1}^m\mathbf{1}_{[-r_i^{\lambda}s_i^{1-\lambda},r_i^{\lambda}s_i^{1-\lambda}]}(x\cdot v_i)dx\ .$$
Set also,
$$G_{\lambda}(\varepsilon):=\int_{x\in\mathbb{R}^n}\int_{u\in\mathbb{R}^m}
\prod_{i=1}^m\mathbf{1}_{[-r_i^{\lambda}s_i^{1-\lambda},r_i^{\lambda}s_i^{1-\lambda}]}(x\cdot v_i+u_i)\varepsilon^{-n}f_{m+n}(u/\varepsilon)dudx\ ,$$
where $\varepsilon>0$, $u_i=u\cdot e_i$, $i=1,\dots, m$ and $\{e_1,\dots,e_m\}$ is an orthonormal basis in $\mathbb{R}^m$.
It follows by the change of variables
$U:=u/\varepsilon$ that
$$G_{\lambda}(\varepsilon):=\int_{x\in\mathbb{R}^n}\int_{U\in\mathbb{R}^m}
\prod_{i=1}^m\mathbf{1}_{[-r_i^{\lambda}s_i^{1-\lambda},r_i^{\lambda}s_i^{1-\lambda}]}(x\cdot v_i+\varepsilon U_i)f_{m+n}(U)dUdx\ $$
$$\xrightarrow{\varepsilon\rightarrow 0^+}  \int_{\mathbb{R}^n}\prod_{i=1}^m\mathbf{1}_{[-r_i^{\lambda}s_i^{1-\lambda},r_i^{\lambda}s_i^{1-\lambda}]}(x\cdot v_i)dx
\int_{U\in \mathbb{R}^m}f_{m+n}(U)dU=\|f_{m+n}|_{\mathbb{R}^m}\|_{1}\cdot|R_{\lambda}| \ .$$
Thus, it suffices to prove that, for $\varepsilon>0$, $$G_{\lambda}(\varepsilon)\geq G_{1}(\varepsilon)^{\lambda}G_{0}(\varepsilon)^{1-\lambda} \ .$$
Using the change of variables $w_i:=\frac{u_i+x\cdot v_i}{s_i}$, $i=1,\dots,m$, we get:
$$G_{\lambda}(\varepsilon)=A\int_{w\in\mathbb{R}^m}\int_{x\in\mathbb{R}^n}
\prod_{i=1}^m\mathbf{1}_{[-r_i^{\lambda}s_i^{1-\lambda},r_i^{\lambda}s_i^{1-\lambda}]}(s_iw_i)f_{m+n}\Big(\varepsilon^{-1}\sum_{i=1}^m( s_iw_i-x\cdot v_i)e_i\Big)dxdw$$
$$=\lim_{a\rightarrow\infty}A\int_{w\in\mathbb{R}^m}\int_{x\in\mathbb{R}^n}
\prod_{i=1}^m\mathbf{1}_{[-r_i^{\lambda}s_i^{1-\lambda},r_i^{\lambda}s_i^{1-\lambda}]}(s_iw_i)\mathbf{1}_{aC_n}(x)
f_{m+n}\Big(\varepsilon^{-1}\sum_{i=1}^m(s_iw_i-x\cdot v_i)e_i\Big)dxdw\ ,
$$
where $A=\varepsilon^{-n}s_1\dots s_n$. Define the (singular) linear map $T:\mathbb{R}^{m+n}\rightarrow \mathbb{R}^{m+n}$, with
$T(w,x)=\varepsilon^{-1}\sum_{i=1}^m(s_iw_i-x\cdot v_i)e_i$. Note that the linear map $T_{\delta}:=T+\delta Id_{\mathbb{R}^{m+n}}$ becomes invertible, for
$\delta>0$, small enough. Therefore, if $$A_a:=\textnormal{diag}\big(\log (r_1s_1^{-1}),\dots ,\log (r_ms_m^{-1}),\log a ,\dots ,\log a\big)\ ,\ a>1\ ,$$
then
\begin{eqnarray*}
G_{\lambda}(\varepsilon)
&=&A\lim_{a\rightarrow \infty}\lim_{\delta\rightarrow 0^+}\int_{e^{[\lambda\cdot 1+(1-\lambda)\cdot 0]A_a}
C_{m+n}}f_{m+n}(T_{\delta}z)dz\\
&=&:A\lim_{a\rightarrow \infty}\lim_{\delta\rightarrow 0^+}F(a,\delta,\lambda)\ .
\end{eqnarray*}
Using our assumption, $$F(a,\delta,\lambda)\geq F(a,\delta,1)^{\lambda}F(a,\delta,0)^{1-\lambda} \ ,$$
for all $a>1$, $\delta>0$ ($\delta$ small enough). This proves our claim. $\Box$
\begin{remark}
\label{remark-section-3}
The case of $f_n$ being the Gaussian density seems to be the most promising in the attempt of proving the log-Brunn-Minkowski inequality.
It follows by the previous Theorem and Theorem \ref{log-BM for Lebesgue implies general log-BM} that the log-Brunn-Minkowski inequality is true in any dimension
and for every log-concave density if and only if it holds true for the (standard) Gaussian density and for parallelepipeds with parallel facets, in
all dimensions.
\end{remark}

\section{The dual log-Brunn-Minkowski inequality}
The main goal of this section is to establish the following dual logarithmic Brunn-Minkowski inequality (see Corollary \ref{dual-lbm-extension-to-quermassintegrals} below).
\begin{theorem}\label{dual-log-BM}
Let $K$, $L$ be two convex bodies in $\mathbb{R}^n$ and $\lambda\in (0,1)$. Then,
\begin{equation}
|(\lambda\cdot K+_0(1-\lambda)\cdot L)^{\circ}|\leq |K^{\circ}|^{\lambda}|L^{\circ}|^{1-\lambda}\ .\label{equation-dual-log-bm}
\end{equation}
\end{theorem}
Once Theorem \ref{dual-log-BM} is established, one can follow Firey's argument \cite{Fi2}
to prove Corollary \ref{dual-lbm-extension-to-quermassintegrals} (see below),
where the volume is replaced by the other quermassintegrals.
Since the dual $L^0$-sum contains the dual $L^p$-sum, for $p\geq 0$, Theorem \ref{dual-log-BM} extends immediately to the $L^p$-setting,
for all $p\geq 0$. Therefore, it is stronger than Firey's \cite{Fi1} dual Brunn-Minkowski inequality.
It is also stronger than the dual Brunn-Minkowski inequality with respect to $L^0$-radial sums, established in \cite{GHWY}.
It seems plausible that the equality cases in (\ref{equation-dual-log-bm}) are exactly the equality cases in Conjecture
\ref{Conjecture-log-bm-lebesque} (here of course non-symmetric bodies are allowed); we do not address this here.

Without loss of generality, we may assume that $K$ and $L$ contain 0 in their interiors. Otherwise, the assertion would be trivial.
As in the previous section, we will prove our claim for the (asymmetric) discrete approximations of the logarithmic sum $K$ and $L$.
The rest of the proof will follow by compactness. Set
$$AR_{\lambda}=\big\{x\in\mathbb{R}^n\ \big|\ x\cdot v_i\leq r_i^{\lambda}s_i^{1-\lambda},\ i=1,\dots,m\big\} \ ,$$
where $v_1,\dots,v_m\in S^{n-1}$, $r_i=h_K(v_i)$, $s_i=h_L(v_i)$, $i=1,\dots,m$. Since $K$ and $L$ contain 0 in their interiors,
it is true that $r_i,s_i>0$,
$i=1,\dots,m$, thus $AR_{\lambda}$ is well defined. One, then, needs to prove that the function $(0,1)\ni\lambda\mapsto
|(AR_{\lambda})^{\circ}|$
is log-convex (i.e. its logarithm is convex). On the other hand,
$$(AR_{\lambda})^{\circ}=\textnormal{conv}\big\{(r_i^{-1})^{\lambda}(s_i^{-1})^{1-\lambda}u_i\ \big|\ i=1,\dots,m\big\} \ ,$$
therefore the proof of Theorem \ref{dual-log-BM} reduces to the proof of the following:
\begin{theorem}\label{dual B-thm}
Let $x_1,\dots,x_m\in \mathbb{R}^n$, $a_1,\dots,a_m\in\mathbb{R}$ and consider the family of polytopes
$$P_t=\textnormal{conv}\big\{e^{a_it}x_i\ \big|\ i=1,\dots,m\big\} \ , \ t\in(t_1,t_2) \ ,$$
for some $t_1<t_2$. If $P_t$ contains the origin in its interior for all $t\in(t_1,t_2)$, then the function $$(t_1,t_2)\ni t\mapsto |P_t|$$
is log-convex.
\end{theorem}
We remark here that Theorem \ref{dual B-thm} may be viewed as the dual version of the B-conjecture for uniform measures in the following sense:
If $K$, $L$ are convex bodies that contain the origin in their interiors, then the function $t\mapsto \big|\big((e^tK)\cap L\big)^{\circ}\big|$
is log-convex.
The idea for the proof of Theorem \ref{dual B-thm} is taken from Saroglou \cite[Theorem 3.1]{Sa2}.
First we will need two easy lemmas.
\begin{lemma}\label{sum of log convex is log convex}
Let $f_1,\dots, f_m:\mathbb{R}^n\rightarrow \mathbb{R}_+$ be log-convex functions. Then, their sum is log-convex.
\end{lemma}
Proof. For $\lambda\in (0,1)$, $t_1,t_2\in \mathbb{R}^n$, we have
$$\sum_{i=1}^mf_i(\lambda t_1+(1-\lambda)t_2)\leq \sum_{i=1}^mf_i^{\lambda}(t_1)f_i^{1-\lambda}(t_2)\leq
\Big[\sum_{i=1}^mf_i(t_1)\Big]^{\lambda}\Big[\sum_{i=1}^mf_i(t_2)\Big]^{1-\lambda} \ .$$
This proves our assertion. $\Box$

\begin{lemma}\label{log-affine}
Let $x_1,\dots,x_n\in \mathbb{R}^n$, $a_1,\dots,a_n\in\mathbb{R}$. Then, the function
$$t\mapsto\big|\textnormal{conv}\big\{0,e^{a_1t}x_1,\dots,e^{a_nt}x_n\big\}\big|$$
is log-affine and therefore log-convex.
\end{lemma}
Proof. We have:
$$\big|\textnormal{conv}\big\{0,e^{a_1t}x_1,\dots,e^{a_nt}x_n\big\}\big|=\frac{1}{n!}\big|\det(e^{a_1t}x_1,\dots,e^{a_nt}x_n)\big|
=\frac{e^{a_1t}\dots e^{a_nt}}{n!}\ ,$$proving our claim. $\Box$
\\
\\
Proof of Theorem \ref{dual B-thm}:\\
We need to prove that the function $(t_1,t_2)\ni t\mapsto  |P_t|$ is log-convex. Actually, we need to prove that for any
$s_1,s_2\in (t_1,t_2)$, $s_1<s_2$, $$\big|P_{\frac{s_1+s_2}{2}}\big|\leq |P_{s_1}|^{1/2}|P_{s_2}|^{1/2}\ .$$
Set $s=(s_1+s_2)/2$, $p=(s_2-s_1)/2$. Let $\{T_1,\dots,T_k\}$ be
a tringulation of the boundary of $P_s$; that is a subdivision of the boundary of $P_s$ into non-overlapping simplices,
whose vertices are vertices of $P_s$.
Set $$\Delta_i:=\textnormal{conv}\big(\{0\}\cup T_i\big)\  , \ i=1,\dots,k \ .$$
Then, the family $\{\Delta_1,\dots,\Delta_k\}$ is a triangulation of $P_s$.
For $i=1,\dots,k$, consider the following transformation of $\Delta_i$:
If
$\Delta_i=\textnormal{conv}\{0,e^{a_{j_1}s}x_{j_1},\dots,e^{a_{j_n}s}x_{j_n}\}$, for some $1\leq j_1<\dots<j_n\leq m$, set
$$\Delta_{i,r}=\textnormal{conv}\{0,e^{a_{j_1}(r+s)}x_{j_1},\dots,e^{a_{j_n}(r+s)}x_{j_n}\}\ , \ \ r\in[-p,p]\ .$$
It is clear that the $\Delta_{i,r}$'s
are non-overalping,
for $r\in[-p,p]$. This is because, for $i=1,\dots,k$, $r\in[-p,p]$, $\Delta_{i,r}$ is contained in the positive cone spanned by
$\Delta_i$
and every two
such cones are, by construction, non-overlapping. Now, it is clear that
\begin{equation}
|P_{s+r}|=\Big|\textnormal{conv}\Big(\bigcup^k_{i=1}\Delta_{i,r}\Big)\Big|\geq\Big(\sum_{i=1}^k|\Delta_{i,r}|\Big) \ ,\ r\in[-p,p] \ .
\label{ineq-of-proof-of Thm-dual-B-theorem}
\end{equation}
By Lemmas \ref{sum of log convex is log convex} and \ref{log-affine}, the function $[-p,p]\ni r \mapsto\sum_{i=1}^k|\Delta_{i,r}|$
is log-convex. Thus,
$$|P_s|=\sum_{i=1}^k|\Delta_{i}|=\sum_{i=1}^k|\Delta_{i,0}|\leq \Big(\sum_{i=1}^k|\Delta_{i,-p}|\Big)^{1/2}\Big(\sum_{i=1}^k|\Delta_{i,p}|\Big)^{1/2}
\leq |P_{s-p}|^{1/2}|P_{s+p}|^{1/2}=|P_{s_1}|^{1/2}|P_{s_2}|^{1/2} ,$$
as required. $\Box$

\begin{corollary}\label{dual-lbm-extension-to-quermassintegrals}
Let $K$, $L$ be two convex bodies that contain 0 in their interior. For $i=1,\dots,n-1$, $p\geq 0$ and $\lambda\in (0,1)$, the following is true:
$$W_i\Big(\big[\lambda\cdot K+_p(1-\lambda)\cdot L\big]^{\circ}\Big)\leq W_i(K^{\circ})^{\lambda}W_i(L^{\circ})^{1-\lambda} \ .$$
\end{corollary}
Corollary \ref{dual-lbm-extension-to-quermassintegrals} also generalizes a result of Firey \cite{Fi2}, who proved this for the $L^1$-sum.
This, was recently extended in the $L^p$-case, for $p\geq 1$ in \cite{HC}, where it was explained that by the homogeneity of
the quermassintegrals,
dual Brunn-Minkowski inequalities have dimension-dependent equivalent forms (in the same manner as the original Brunn-Minkowski
inequality does; see e.g. \cite{Ga1}).
Adoupting the same argument we obtain:
\begin{corollary}\label{dimension-dependent}
Let $K$, $L$ be two convex bodies that contain the origin in their interior. If $i\in\{1,\dots,n-1\}$, $p\geq 0$ and $\lambda\in (0,1)$, then
$$W_i\Big(\big[\lambda\cdot K+_p(1-\lambda)\cdot L\big]^{\circ}\Big)^{-\frac{p}{n-i}}\leq\lambda W_i(K^{\circ})^{-\frac{p}{n-i}}+(1-\lambda) W_i(L^{\circ})^{-\frac{p}{n-i}}\ .$$
\end{corollary}
Proof. Use Corollary \ref{dual-lbm-extension-to-quermassintegrals} with
$\overline{K}=W_i(K^{\circ})^{1/n-i}$, $\overline{L}=W_i(L^{\circ})^{1/n-i}$,
$\overline{\lambda}=\lambda W_i(L^{\circ})^{-p/n-i}\cdot\mu^{-1}$,
where $\mu=\lambda W_i(K^{\circ})^{-\frac{p}{n-i}}+(1-\lambda) W_i(L^{\circ})^{-\frac{p}{n-i}}$,
in the place of $K$, $L$, $\lambda$ respectively. $\Box$\\
\\
For the proof of Corollary \ref{dual-lbm-extension-to-quermassintegrals}, the following (contained in an earlier version of \cite{Sa1}) is required.
\begin{lemma}\label{section contained}
Let $\lambda \in [0,1]$, $K$, $L$ be convex bodies in $\mathbb{R}^n$ and $H$ be a subspace of $\mathbb{R}^n$.
Then, $$\lambda\cdot(K\cap H)+_0(1-\lambda)\cdot(L\cap H)\subseteq [\lambda\cdot K+_0(1-\lambda)\cdot L]\cap H \ ,$$
where the logarithmic sum in the first part of the previous inclusion is considered with respect to the subspace $H$.
\end{lemma}
Proof. Note that if $x,\ u\in H$ and $y\in H^{\bot}$, then $x\cdot (u+y)=x\cdot u$, $$
h_{K\cap H}(u+y)=\max_{z\in K\cap H}z\cdot (u+y)=\max_{z\in K\cap H}z\cdot u=h_{K\cap H}(u)$$
and, similarly, $h_{L\cap H}(u+y)=h_{L\cap H}(u)$. Thus, $\lambda\cdot (K\cap H)+_0(1-\lambda)\cdot(L\cap H)$
\begin{eqnarray*}
&=&\{x\in H \ | \ x\cdot(u+y)\leq h^{\lambda}_{K\cap H}(u+y)h^{1-\lambda}_{L\cap H}(u+y),\textnormal{ for all }u\in H, \ y\in H^{\bot}\}\\
&=& \{x\in H \ | \ x\cdot w\leq h^{\lambda}_{K\cap H}(w)h^{1-\lambda}_{L\cap H}(w),\textnormal{ for all }w\in \mathbb{R}^n\}\\
&\subseteq& \{x\in H \ | \ x\cdot w\leq h^{\lambda}_{K}(w)h^{1-\lambda}_{L}(w),\textnormal{ for all }w\in \mathbb{R}^n\}\\
&=& [\lambda\cdot K+_0(1-\lambda)\cdot L]\cap H \ . \ \ \ \Box
\end{eqnarray*}
Proof of Corollary \ref{dual-lbm-extension-to-quermassintegrals}:\\
We will make use of Firey's argument for passing from the volume to other quermassintegrals (see \cite{Fi2}) and
the fact that the $p$-convex combination of convex bodies contains the logarithmic convex combination, $p>0$.
It follows immediately by Lemma \ref{section contained}, that
$$\big[\lambda\cdot K+_0(1-\lambda)\cdot L\big]^{o}\big|H\subseteq \big[\lambda\cdot(K\cap H)+_0(1-\lambda)\cdot(L\cap H)\big]^{\circ}\ .$$
Therefore, by the Kubota formula and Theorem \ref{dual-log-BM}, we obtain:
\begin{eqnarray*}
W_i\Big(\big[\lambda\cdot K+_p(1-\lambda)\cdot L\big]^{\circ}\Big)&\leq&W_i\Big(\big[\lambda\cdot K+_o(1-\lambda)\cdot L\big]^{\circ}\Big)\\
&=&\frac{\omega_n}{\omega_{n-i}}
\int_{{\cal{G}}_{n,n-i}}\Big|\big[\lambda\cdot K+_0(1-\lambda)\cdot L\big]^{o}\big|H\Big|_{n-i}dH \\
&\leq&\frac{\omega_n}{\omega_{n-i}}\int_{{\cal{G}}_{n,n-i}}\Big|\big[\lambda\cdot(K\cap H)+_0(1-\lambda)\cdot(L\cap H)\big]^{\circ}\Big|_{n-i}dH\\
&\leq&\frac{\omega_n}{\omega_{n-i}}\int_{{\cal{G}}_{n,n-i}}|(K\cap H)^{\circ}|_{n-i}^{\lambda}|(L\cap H)^{\circ}|_{n-i}^{1-\lambda}dH\\
&=&\frac{\omega_n}{\omega_{n-i}}\int_{{\cal{G}}_{n,n-i}}|K^{\circ}|H|_{n-i}^{\lambda}|L^{\circ}|H|_{n-i}^{1-\lambda}dH\\
&\leq&\bigg[\frac{\omega_n}{\omega_{n-i}}\int_{{\cal{G}}_{n,n-i}}|K^{\circ}|H|_{n-i}dH\bigg]^{\lambda}
\bigg[\frac{\omega_n}{\omega_{n-i}}\int_{{\cal{G}}_{n,n-i}}|L^{\circ}|H|_{n-i}dH\bigg]^{1-\lambda}\\
&=&W_i(K^{\circ})^{\lambda}W_i(L^{\circ})^{1-\lambda} \ ,\ i=1,\dots,n-1 \ .
 \ \ \ \ \Box\end{eqnarray*}\\
\\
Before ending this note, we would like to state a consequence of Theorem \ref{dual-log-BM} that concerns the logarithmic sum itself,
rather than its dual.

\begin{corollary}\label{dual-last-corollary}
Let $\Delta_1,\Delta_2\subseteq \mathbb{R}^2$ be two triangles whose centroids are at the origin. Then,
\begin{equation}
|\lambda\cdot \Delta_1+_0(1-\lambda)\cdot\Delta_2|\geq |\Delta_1|^{\lambda}|\Delta_2|^{1-\lambda} \ .
\label{eq-triangles}
\end{equation}
\end{corollary}
Proof. It is well known (see \cite{Ma}) that if $K$ is any planar convex body, then
\begin{equation}
|K|\cdot|K^{\circ}|\geq |\Delta_1|\cdot |\Delta_1^{\circ}|\ ,\label{volume-product-planar}
\end{equation}
Now, if (\ref{eq-triangles}) is not true, then by Theorem \ref{dual-log-BM} 
and (\ref{volume-product-planar}) we get:
$$|\lambda\cdot \Delta_1+_0(1-\lambda)\cdot\Delta_2|\cdot|(\lambda\cdot \Delta_1+_0(1-\lambda)\cdot\Delta_2)^{\circ}|<|\Delta_1|^{\lambda}|\Delta_2|^{1-\lambda}
|\Delta_1^{\circ}|^{\lambda}|\Delta_2^{\circ}|^{1-\lambda}=|\Delta_1|\cdot |\Delta_1^{\circ}| \ ,$$
which contradicts (\ref{volume-product-planar}). $\Box$
\\
\\
{\bf Agnowledgement.} I would like to thank A. Marsiglietti for reading carefully this manuscript and for many useful commends,
especially for pointing me Corollary \ref{Gaussian-BM}. I would also like to thank A. Zvavitch and O. Guedon for explaining 
me some very nice things.

\bibliographystyle{plain}

\bigskip

\noindent \textsc{Ch.\ Saroglou}: Department of Mathematics,
Texas A$\&$M University, 77840 College Station, TX, USA.

\smallskip

\noindent \textit{E-mail:} \texttt{saroglou@math.tamu.edu \ \&\ christos.saroglou@gmail.com}

\end{document}